\documentclass[10pt]{amsart}
\usepackage{pstricks,pstcol,pst-plot}
\parskip=\smallskipamount
\numberwithin{equation}{section}

\definecolor{myblue}{rgb}{0.66,0.78,1.00}
\definecolor{Apricot}{cmyk}{0,0.25,1,0}             
\definecolor{OrangeRed}{cmyk}{0,0.6,1,0}            
\definecolor{Lavender}{cmyk}{0,0.48,0,0}
\definecolor{Violet}{cmyk}{0.79,0.88,0,0}
\definecolor{DarkBlue}{cmyk}{1,1,0,0.20}
\definecolor{Black}{cmyk}{0,0,0,1}

\newtheorem{theorem}{Theorem}[section]
\newtheorem{lemma}[theorem]{Lemma}

\theoremstyle{definition}

\newtheorem{remark}[theorem]{Remark}

\newcommand{\C}{\mathbb{C}}

\renewcommand{\P}{\mathbb{P}}
\newcommand{\Z}{\mathbb{Z}}
\newcommand{\R}{\mathbb{R}}

\newcommand{\cC}{\mathcal{C}}
\newcommand{\cH}{\mathcal{H}}
\newcommand{\cL}{\mathcal{L}}

\newcommand\wt{\widetilde}
\newcommand{\bs}{\backslash}
\newcommand{\e}{\epsilon}

\begin{document}
\title[Deformations  of Stein structures]
{Deformations of Stein structures \\ and extensions of holomorphic mappings}
\author{Franc Forstneri\v c and Marko Slapar}
\address{Institute of Mathematics, Physics and Mechanics \\ University of Ljubljana \\ Jadranska 19, 1000 Ljubljana, Slovenia}
\email{franc.forstneric@fmf.uni-lj.si, marko.slapar@fmf.uni-lj.si }

\thanks{Supported by grants P1-0291 and J1-6173, Republic of Slovenia.}

%
%
\subjclass[2000]{32H02, 32Q28, 32Q30, 32Q55, 32Q60, 32T15, 57R17}
\date{Sept.\ 19, 2005} 
\keywords{Stein manifolds, complex structures, holomorphic mappings}

\begin{abstract}
Assume that $A$ is a closed complex subvariety of a Stein manifold $X$
and that $f\colon X\to Y$ is a continuous map to a complex manifold $Y$
such that the restriction $f|_A\colon A\to Y$ is holomorphic
on $A$. After a homotopic deformation of the Stein structure 
outside a neighborhood of $A$ in $X$ we find a holomorphic map 
$\wt f\colon X\to Y$ which agrees with $f$ on $A$ and which 
is homotopic to $f$ relative to $A$. When $\dim_\C X=2$ we must also 
change the $\cC^\infty$ structure on $X\bs A$. 
\end{abstract}
\maketitle

\section{Introduction} 
According to a classical theorem of H.\ Cartan 
every holomorphic function on a closed complex subvariety 
of a Stein manifold $X$ extends to a holomorphic function on 
all of $X$ (\cite{GR}, \cite{Ho}). This extension property 
fails in general for mappings $X\to Y$ to more general 
complex manifolds unless $Y$ enjoys a certain holomorphic 
flexibility property introduced in \cite{FCAP} and \cite{FFourier2}.
In this paper we show that Cartan's extension theorem  holds
for maps to an arbitrary complex manifold if we allow 
homotopic deformations of the complex structure
(and of the underlying smooth structure when $\dim_\C X=2$)
in the complement of the given subvariety in the source Stein manifold $X$. 
The following is a simplified version of theorem 
\ref{T3.1} in \S 3 below.

%
%
%
%
\begin{theorem}
\label{Tmain1}
Let $X$ be a Stein manifold with $\dim_\C X\ne 2$ 
and let $A$ be a closed complex subvariety of $X$. Given a 
continuous map $f\colon X\to Y$ to a complex manifold $Y$ 
such that $f|_A\colon A\to Y$ is holomorphic, 
there is a homotopy $(J_t,f_t)_{t\in [0,1]}$, 
consisting of integrable complex structures $J_t$ on $X$
and of continuous maps $f_t\colon X\to Y$, satisfying
the following properties:
\begin{itemize}
\item[(i)] $J_0$ is the initial complex structure on $X$, 
$J_t=J_0$ in a neighbor\-hood of $A$ for each $t\in[0,1]$, 
and $J_1$ is a Stein structure on $X$;
\item[(ii)] $f_0=f$, $f_t|_A=f|_A$ for every $t\in [0,1]$,
and $f_1\colon X\to Y$ is $J_1$-holomorphic.
\end{itemize}
\end{theorem}

This result is a relative version (with interpolation)
of theorem 1.1 in \cite{FS} to the effect  that every continuous map
$f\colon X\to Y$ from a Stein manifold $(X,J)$ of complex dimension 
$\ne 2$ is homotopic to a map 
$\wt f\colon X\to Y$ which is holomorphic with respect to 
some Stein structure $\wt J$ on $X$ that is homotopic to $J$
through a family of integrable (but not necessarily Stein) complex
structures on $X$. 
When $\dim_\C X=2$ it is in general also necessary to change 
the smooth structure on $X$; see theorem \ref{T4.1}.

The first author proved in \cite{FFourier2} that
the conclusion of theorem \ref{Tmain1} holds for all data 
$(X,A,f)$ without changing the Stein structure 
on $X$ if and only if the target manifold $Y$ satisfies
the {\em convex approximation property} (CAP), introduced
in \cite{FCAP}, to the effect that every holomorphic map
$K\to Y$ from a compact convex set $K\subset \C^n$,
$n=\dim X+\dim Y$, is a 
uniform limit of entire maps $\C^n\to Y$. (See also \cite{La}.)
Among the conditions implying CAP we mention
complex homogeneity and, more generally, the existence of a 
finite dominating family of holomorphic sprays. For a 
more complete discussion of this subject see \cite{Fflex}.

A Stein structure $J_1$ in theorem \ref{Tmain1} 
can be chosen such that $(X,J_1)$ is biholomorphic to $(\Omega,J|_{T\Omega})$
for some $J$-Stein domain  $\Omega\subset X$ which contains $A$ 
and is diffeotopic to $X$ relative to $A$. 
Here is the precise result; for $A=\emptyset$
this is theorem 1.2 in \cite{FS}.

\begin{theorem}
\label{Tmain2}
Let $X$ be a Stein manifold with $\dim_\C X\ne 2$ and 
let $A$ be a closed complex subvariety of $X$.
Given a  continuous map $f\colon X\to Y$ to a complex manifold $Y$ 
such that $f|_A\colon A\to Y$ is holomorphic, there exist a 
Stein domain $\Omega\subset X$ containing $A$, a holomorphic map 
$\wt f\colon \Omega\to Y$, and a diffeomorphism $h\colon X\to \Omega$
which is diffeotopic to $id_X$ by a diffeotopy that is
fixed on a neighborhood $A$, such that the map $\wt f\circ h\colon X\to Y$
is homotopic to $f$ relative to $A$.
\end{theorem}

Theorem \ref{Tmain1} is implied by theorem \ref{Tmain2} as follows.
Let $h_t\colon X \to h_t(X)\subset X$ be a diffeotopy  as
in theorem \ref{Tmain2}, satisfying $h_0=id_X$, $h_1=h \colon X\to\Omega$, 
and such that $h_t$ is the identity map in a neighborhood of $A$ 
for each $t\in[0,1]$. Let $J_t=h_t^*(J)$ denote the (unique) complex structure 
on $X$ satisfying $dh_t\circ J_t = J\circ dh_t$ 
on $TX$. The homotopy $\{J_t\}_{t\in [0,1]}$ then satisfies 
theorem \ref{Tmain1} (i), and the map $f_1:=\wt f\circ h\colon X\to Y$ 
is $J_1$-holomorphic and satisfies theorem \ref{Tmain1} (ii).

\begin{remark}
Although a Stein structure $J_1$ satisfying the conclusion of theorem \ref{Tmain1}
must in general depend on the initial map $f$, 
we can choose the same $J_1$ for all members of a compact Hausdorff family of maps;
this can be seen by applying the parametric versions of the main tools as in \cite{FS}.
The analogous remark applies to theorem \ref{Tmain2} in which 
the Stein domain $\Omega \subset X$ can be chosen
the same for all maps in a compact Hausdorff family. 
\end{remark}

Theorems \ref{Tmain1} and \ref{Tmain2} are proved in \S 3. 
The main inductive step is furnished by lemma \ref{Mainlemma}
whose proof relies on the tools developed by Eliashberg \cite{E} 
and the authors \cite{FS}. The underlying geometric construction 
in this paper is more intricate than the one in \cite{FS} 
due to the presence of a subvariety.  
In \S 4 we discuss the analogous result for maps from Stein surfaces
($\dim_\C X=2$), using results of Gompf \cite{Go1}, \cite{Go2}.

%
%
%
%
\section{The main lemma}
An {\em almost complex structure} on an even dimensional 
smooth manifold $X$ is a smooth endomorphism
$J\in \mathrm{End}_\R(TM)$ satisfying $J^2=-Id$.
It gives rise to the conjugate differential $d^c$, defined on 
functions by $\langle d^c\rho,v\rangle= - \langle d\rho, Jv\rangle$
for $v\in TX$ (equivalently, $d^c=-J^*d$), and the Levi form operator $dd^c$.
The structure $J$ is {\em integrable} if every point of $X$ 
admits an open neighborhood $U\subset X$ and a 
$J$-holomorphic coordinate map of maximal rank
$z=(z_1,\ldots,z_n) \colon U  \to \C^n$ $(n=\frac{1}{2}\dim_\R X$),
i.e., satisfying $dz\circ J=idz$; for a necessary and sufficient 
integrability condition see \cite{NN}. 

Let $(X,J_X)$ and $(Y,J_Y)$ be a pair of (almost) complex manifolds.
A smooth map $f\colon X\to Y$ is $(J_X,J_Y)$-holomorphic 
if $df\circ J_X = J_Y\circ df$. Since the complex structure
on $Y$ will be kept fixed in our proofs,
we shall simply say that $f$ is $J_X$-holomorphic.

We assume familiarity with standard complex analytic notions 
such as (strong) plurisubharmonicity and (strong) pseudoconvexity 
(see \cite{GR}, \cite{Ho}). Since we shall deal with several
different complex structures on the same manifold, 
we shall often write $J$-holo\-mor\-phic, 
$J$-Stein, $J$-plurisubharmonic, $J$-pseudoconvex, 
etc.

If $(X,J)$ is a Stein manifold and 
$K\subset L\subset X$, with $K$ compact, we shall say
that $K$ is {\em $J$-holomorphically convex in $L$} if 
for every $p\in L\bs K$ there is a $J$-holomorphic 
function $f$ on an open set in $X$ containing $L$ such that 
$|f(p)|>\sup_{x\in K} |f(x)|$. When this holds
with $L=X$, we say that $K$ is $\cH(X,J)$-convex.

The following lemma is the main step
in the proof of theorems \ref{Tmain1} and \ref{Tmain2}.

\begin{lemma}
\label{Mainlemma}
Let $(X,J_X)$ be a Stein manifold with $\dim_\C X=n\ne 2$.
Let $\rho\colon X\to \R$ be a smooth strongly $J_X$-plurisubharmonic 
Morse exhaustion function, and let $c'<c$ be regular values of $\rho$.
Set $K=\{x\in X\colon\rho(x)\le c'\}$, $L=\{x\in X\colon\rho(x)\le c\}$. 
Let $A$ be a closed complex subvariety of $X$.

Assume that $J$ is an almost complex structure on $X$
which is integrable in an open neighborhood $U\subset X$ of $A\cup K$
such that $J$ agrees with $J_X$ in a neighborhood of $A$ in $X$ and 
$K$ is a strongly $J$-pseudoconvex domain with $J$-Stein interior. 

Let $Y$ be a complex manifold endowed with a distance function $d_Y$ 
induced by a Riemannian metric.
Given a continuous map $f\colon X\to Y$ which is $J$-holomorphic 
in a neighborhood of $K$ and such that $f|_A\colon A\to Y$ 
is $J_X$-holomorphic, there exists for every $\e>0$ 
a homotopy of pairs $(J_t,f_t)$ $(t\in[0,1])$, where $J_t$ is
an almost complex structure on $X$ and $f_t\colon X\to Y$
is a continuous map, satisfying the following:
\begin{itemize}
\item[(i)] $J_t$ agrees with $J_0=J$ in a neighborhood of $A\cup K$ 
for all $t\in [0,1]$,
\item[(ii)]  $J_1$ is integrable in a neighborhood of $A\cup L$,
\item[(iii)]  $L$ is a strongly $J_1$-pseudoconvex domain with 
$J_1$-Stein interior, and the set $K$ is $J_1$-holomorphically convex in $L$,  
\item[(iv)] $f_0=f$ and $f_t|_A=f|_A$ for each $t\in[0,1]$,
\item[(v)]  for every $t\in [0,1]$ the map $f_t$ is $J$-holomorphic 
in a neighborhood of $K$ and satisfies
$\sup_{x\in K} d_Y\bigl(f_t(x),f(x)\bigr) <\e$, and
\item[(vi)]   the map $f_1\colon X\to Y$ is $J_1$-holomorphic in
a neighborhood of $L$.
\end{itemize}
If $J$ is integrable on $X$ then all structures $J_t$ $(t\in[0,1])$
can be chosen integrable.
\end{lemma}

\begin{figure}[ht]
\psset{unit=0.6cm} 
\begin{pspicture}(-7,-4)(7,4)

\pscircle[linecolor=OrangeRed,fillstyle=crosshatch,hatchcolor=yellow](0,0){4}   
\rput(0,3){$L$}						                        

\pscircle[linecolor=DarkBlue,fillstyle=crosshatch,hatchcolor=myblue](0,0){2}    
\rput(0,-0.5){$K$}					         	        

\pscurve[linecolor=red,linewidth=1.5pt](-6.5,0)(-4,0.5)(-3,-0,3)(-2,0)(-1,0.3)(0,0.5)
(1,0.3)(2,-0.2)(3,-0.1)(4,0.2)(6.5,0)                                           

%
%
\psarc[linestyle=dashed,linewidth=0.5pt,linecolor=Black](0,0){2.3}{13}{170}                                      
\psecurve[linestyle=dashed,linewidth=0.5pt](2.1,1)(2.2,0.6)(2.8,0.2)(3.4,0.35)(4,0.45)(5,0.5)(6.5,0.25)(7,0.4)        
\psecurve[linestyle=dashed,linewidth=0.5pt](-2,2)(-2.27,0.4)(-2.5,0.2)(-3,0.4)(-4,0.8)(-5,0.8)(-6.5,0.3)(-7,0)   
\psarc[linestyle=dashed,linewidth=0.5pt](0,0){2.3}{195}{340}                                        
\psecurve[linestyle=dashed,linewidth=0.5pt](-2,-1)(-2.2,-0.6)(-2.5,-0.3)(-3,-0.3)(-4,0.2)(-5,0)(-6.5,-0.3)(-7,-0.6)  
\psecurve[linestyle=dashed,linewidth=0.5pt](2,-1.5)(2.18,-0.8)(2.6,-0.4)(3,-0.4)(4,-0.2)(5,-0.2)(6.5,-0.4)(7,-0.4)   

%
%
\rput(-5.5,1.5){$A$}
\psline[linewidth=0.2pt]{->}(-5.5,1.2)(-5.5,0.35)
\rput(5.1,1.5){$U$}
\psline[linewidth=0.2pt]{->}(5,1.2)(5,0.5)

\end{pspicture}
\caption{The main lemma.} 
\label{Fig1}
\end{figure}
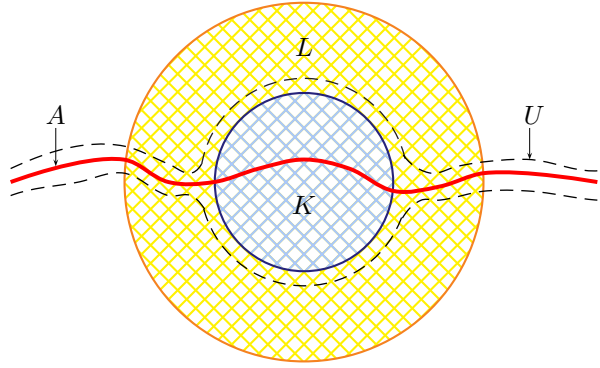

The situation is illustrated on fig.\ \ref{Fig1}:
$J$ is integrable in $U\supset A\cup K$ (shown with the dashed line), 
$f|_A$ is holomorphic with respect to the complex structure on $A$
induced by $J_X$, and $f$ is $J$-holomorphic in a neighborhood of $K$. 
The pair final $(J_1,f_1)$  enjoys the analogous properties on the 
larger set $A\cup L$.

%
%
%
%

\begin{proof}
We may assume that $K=\{\rho\le -1\}$ and $L=\{\rho\le 0\}$. 

Our first goal is to replace $f$ by another map $X\to Y$
which is holomorphic in an open neighborhood of $A\cup K$ in $X$,
without changing $f$ on $A$ (where it is already holomorphic)
and changing it arbitrary little on $K$.

The set $K$, being strongly $J$-pseudoconvex with $J$-Stein interior,
admits a basis of $J$-Stein neighborhoods.
Also, since $K$ is $J_X$-holomorphically convex in $X$
and $J=J_X$ in a neighborhood of $A$, it follows that
$A\cap K$ is holomorphically convex in $A$ with respect
to the complex structure induced by $J$ (or by $J_X$
since they agree near $A$).
Theorem 2.1 in \cite{FFourier2}, applied to the set 
$A\cup K$ in the complex manifold $(U,J|_{TU})$, shows that
$A\cup K$ admits a fundamental basis of open $J$-Stein 
neighborhoods $V_j \subset U$ such that 
$K$ is $J$-holomorphically convex in $V_j$. Replacing $U$
by such a neighborhood we shall assume that $U$ is $J$-Stein
and $K$ is $\cH(U,J)$-convex.

Theorem 3.1 in \cite{FFourier2} now furnishes an open neighborhood 
$U'\subset U$ of $A\cup K$ and a $J$-holomorphic map $f'\colon U'\to Y$ 
such that $f'|_{A}=f|_{A}$ and $f'|_K$ is uniformly as close as desired to $f|_K$.
If the approximation is sufficiently close then it is possible to 
patch $f'$ and $f$ outside a small open neighborhood of $A\cup K$ 
and thereby extend $f'$ to all of $X$. The change from $f$ to $f'$
is accomplished by a homotopy with the desired properties;
replacing $f$ by $f'$ and shrinking $U$ 
we may therefore assume $f\colon X\to Y$ is 
$J$-holomorphic in a $J$-Stein domain $U \supset A\cup K$.

Let $g_1,\ldots,g_r$ be $J_X$-holomorphic functions on $X$ such that 
$$
	A=\{x\in X\colon g_1(x)=0,\ldots,g_r(x)=0\}.
$$
We may assume that $\sum_{j=1}^r |g_j|^2 <1$ on $K$.
For every $\delta>0$ the function  
$$ 
	\phi_\delta = (\rho + 1) + \delta\cdotp \log\biggl( \sum_{j=1}^r |g_j|^2 \biggr)
$$
is strongly $J_X$-plurisubharmonic on $X$, $\phi_\delta<0$ on $K$,
and $A=\{\phi_\delta=-\infty\}$.
A generic choice of $\delta$ insures that 
$\Sigma_\delta := \{x\in L\colon \phi_\delta(x)=0\}$ is a smooth
strongly $J_X$-pseudoconvex hypersurface intersecting $bL$ transversely.

We wish to smoothen the corner of the set 
$\{x\in L\colon \phi_\delta(x)\le 0\}$ along $\Sigma_\delta \cap bL$
so that the new domain will have $J$-Stein interior
and smooth strongly $J$-pseudoconvex boundary. Let 
$\tau_\delta =rmax(\rho,\phi_\delta)$,
where $rmax$ denotes a regularized maximum function (see lemma 5 in \cite{De}). 
The function $\tau_\delta$ is smooth and strongly $J_X$-pluri\-subharmonic
on $X$ (since $rmax$ preserves this property), 
it equals $\rho$ near  $A$ (since $\phi_\delta|_A=-\infty$), 
and it equals $\phi_\delta$ on $\{x \in L\colon \phi_\delta\ge 0\}$ 
(since $\rho\le 0$ on $L$). 
The set $E_\delta= \{x\in L\colon \tau_\delta(x) \le 0\}$
has smooth strongly $J_X$-pseudoconvex boundary 
which coincides with $bL$ in a neighborhood of $A\cap bL$,
and it coincides with $\Sigma_\delta$ in $\{\rho\le c\}$
for some $c<0$ close to $0$. (The set $E_\delta$ is shown as $D_0$
in fig.\ \ref{Fig2} below.) We have $K\subset {\rm Int}E_\delta$ for every $\delta>0$.
As $\delta$ decreases to $0$, $E_\delta$ shrinks down to $K\cup (A\cap L)$. 

We claim that for a sufficiently small $\delta>0$ the set
$E_\delta$ has $J$-Stein interior and strongly $J$-pseudoconvex 
boundary $bE_\delta$. Since $E_\delta$ is contained in the $J$-Stein manifold $U$,
it suffices to verify the latter property; the first one will then
follow from the general theory.  Recall that $J=J_X$ in an open 
set $V\supset A$. The part of $bE_\delta$ which belongs 
to $V$ is strongly $J$-pseudoconvex since $J=J_X$ in $V$.  
The remaining part $bE_\delta \cap (L\bs V)$ converges to $bK\bs V$ in 
the $\cC^\infty$ topology as $\delta$ decreases to $0$
as is seen from the definition of $\phi_\delta$.
Since $bK$ is assumed strongly $J$-pseudoconvex,  
$bE_\delta\bs V$ is also such provided that 
$\delta>0$ is chosen sufficiently small. This establishes the claim.

We now fix a $\delta>0$ satisfying  the above
requirements and set $\tau=\tau_\delta$, $E=E_\delta$. 
We proceed as in the proof of Lemma 6.9 in \cite{ACTA}.
For $t\in [0,1]$ we set 
$$
	\rho_t= (1-t)\tau + t\rho, \quad D_t=\{x\in X\colon \rho_t(x)\le 0\}.
$$
The function $\rho_t$, being a convex combination of two strongly 
$J_X$-pluri\-sub\-harmonic functions $\rho_0=\tau$ and $\rho_1=\rho$, 
is itself strongly $J_X$-plurisubhar\-monic. The sets $D_t$ are strongly 
$J_X$-pseudoconvex with smooth boundaries, except at points in $bD_t$ 
where $d\rho_t=0$. 
We have $E=D_0\subset D_t\subset D_1=L$ for every $t\in [0,1]$; 
as $t$ increases from $0$ to $1$, the domains $D_t$ monotonically increase 
from $D_0$ to $D_1=L$ (fig.\ \ref{Fig2}).

%
%
%
%
\begin{figure}[ht]
\psset{unit=0.7cm} 
\begin{pspicture}(-8,-5)(8,5)


\pscircle[linecolor=OrangeRed,fillstyle=crosshatch,hatchcolor=yellow](0,0){5}       
\rput(0,4.3){$L=D_1$}						                            

\pscircle[linecolor=DarkBlue,fillstyle=crosshatch,hatchcolor=myblue](0,0){2.6}      
\rput(0,0.5){$K$}					         	            

\psline[linecolor=red,linewidth=1.5pt](-8,0)(8,0)                                   
\rput(-7,0.3){$A$}
\rput(7,0.3){$A$}

%
%
\psarc[linewidth=1.2pt](0,0){5}{-3}{3}
\psarc[linewidth=1.2pt](0,0){5}{177}{183}

%
%
\psarc[linewidth=1.2pt](0,0){2.8}{20}{160}
\psarc[linewidth=1.2pt](0,0){2.8}{200}{340}

%
%
%
\psecurve[linewidth=1.2pt](3,4)(2.6,1)(3,0.7)(4,0.7)(4.97,0.25)(4,-1)
\psecurve[linewidth=1.2pt](3,-4)(2.6,-1)(3,-0.7)(4,-0.7)(4.97,-0.25)(4,1)
\psecurve[linewidth=1.2pt](-3,4)(-2.6,1)(-3,0.7)(-4,0.7)(-4.97,0.25)(4,-1)
\psecurve[linewidth=1.2pt](-3,-4)(-2.6,-1)(-3,-0.7)(-4,-0.7)(-4.97,-0.25)(-4,1)

\rput(-3.8,0.3){$D_0$}   
\rput(3.8,0.3){$D_0$}   

%
%
%
\psarc[linestyle=dashed,linewidth=0.8pt](0,0){3.6}{30}{150}
\psarc[linestyle=dashed,linewidth=0.8pt](0,0){3.6}{210}{330}

%
%
%
\psecurve[linestyle=dashed,linewidth=0.8pt](2,3)(3.15,1.75)(3.6,1.4)(4.4,1)(5,0.25)(4,-1)
\psecurve[linestyle=dashed,linewidth=0.8pt](-2,3)(-3.15,1.75)(-3.6,1.4)(-4.4,1)(-5,0.25)(-4,-1)
\psecurve[linestyle=dashed,linewidth=0.8pt](2,-3)(3.15,-1.75)(3.6,-1.4)(4.4,-1)(5,-0.25)(4,1)
\psecurve[linestyle=dashed,linewidth=0.8pt](-2,-3)(-3.15,-1.75)(-3.6,-1.4)(-4.4,-1)(-5,-0.25)(-4,1)

\rput(0,3.2){$D_t$}

\psline[linewidth=0.2pt]{<-}(3.05,2.05)(4.4,3.4)
\rput(5.6,3.8){$\{\rho_t=0\}=bD_t$}

\end{pspicture}
\caption{The sets $D_t=\{\rho_t\le 0\}$.} 
\label{Fig2}
\end{figure}
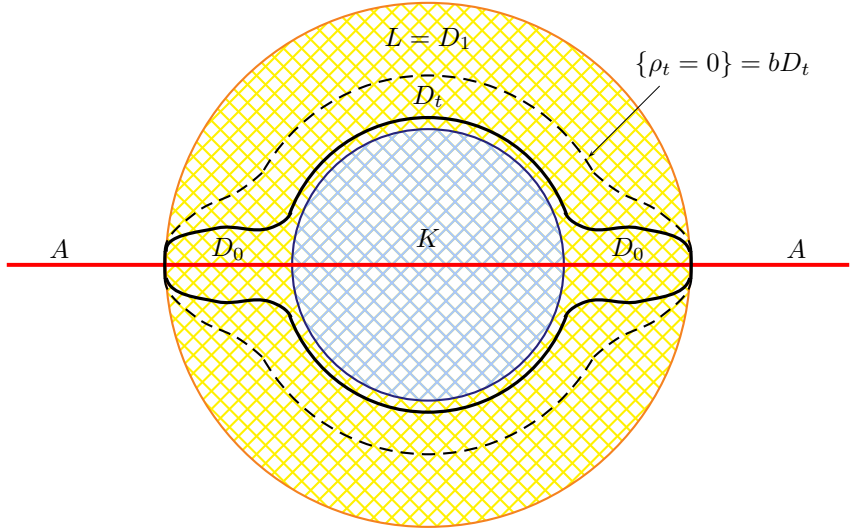

Our goal is to show that the domain $L=D_1$ can be obtained
(up to a diffeomorphism) from the domain $D_0$ by attaching handles of 
indices $\le n$. To this end we investigate the singular points 
of the hypersurfaces $bD_t =\{\rho_t=0\}$ for $t\in[0,1]$. By the construction,
all these boundaries coincide on $\{\rho=0,\tau = 0\}=bL\cap bD_0$,
and this set is a relative neighborhood of $A\cap bL$ in $bL$.
Since the boundaries $bD_0=\{\tau= 0\}$ and $bL=\{\rho = 0\}$ are 
smooth and they intersect transversely along $bL\cap bD_0$, we
see that all nonsmooth points of $bD_t$ are contained in the open set
$
	\Omega= \{\rho<0,\ \tau>0\}= {\rm Int}L \bs D_0.
$
The defining equation of $D_t \cap \Omega$ is 
$\tau \le t(\tau-\rho)$; dividing by $\tau-\rho>0$ we get 
$$   
	D_t\cap \Omega = \{x\in \Omega \colon 
	\sigma(x) := 	\frac{\tau(x)}{\tau(x)-\rho(x)} \le  t\},\quad  t\in[0,1].
$$
The critical point equation $d\sigma =0$ is equivalent to 
$$
	(\tau -\rho)d\tau - \tau(d\tau -d\rho)= \tau d\rho-\rho d\tau =0.
$$
Generic choices of $\rho$ and $\tau$ insure that there are
at most finitely many solutions $p_1,\ldots,p_m\in \Omega$,
all nondegenerate (Morse) and belonging to pairwise distinct level sets 
of $\sigma$, and there are no solution on $b\Omega$. 
A calculation gives the following relationship between the 
$J_X$-Levi forms of these functions at a critical point $p_j$ of $\sigma$:  
$$
	\bigl(\tau(p_j) -\rho(p_j) \bigr)^2 \cL_\sigma(p_j) 
	= \tau(p_j) \cL_\rho(p_j) - \rho(p_j) \cL_\tau(p_j).
$$
(In local holomorphic coordinates $z=(z_1,\ldots,z_n)$ at $p$,
and with $w\in\C^n$, we have $\cL_\sigma(p)\cdotp w = \sum_{j,k=1}^n 
\frac{\partial\sigma}{\partial z_j\partial \overline z_k}(p) w_j\overline w_k$,
and $\cL_\sigma(p)>0$ means that this expression is positive for every
$w\ne 0$.) Since $\tau(p_j)>0$, $-\rho(p_j)>0$ and the functions
$\tau$ and $\sigma$ are strongly $J_X$-plurisubharmonic, we obtain $\cL_\sigma(p_j)>0$. 
It follows that the Morse index of $\sigma$ at $p_j$ is $\le n=\dim_\C X$. 
(If not, the $\R$-linear subspace $\Lambda$ of $T_{p_j} X$,
corresponding to all the negative eigenvalues of the real Hessian of 
$\sigma$ at $p_j$, would have real dimension at least $n+1$
and hence $\Lambda\cap J_X(\Lambda)$ would contain a 
complex line $\lambda$; the restriction of 
$\cL_\sigma(p_j)$ to $\lambda$ would therefore 
be negative, a contradiction.)  

Choose numbers $t_0=0<t_1<t_2<\ldots <t_m=1$ which are regular values
of $\sigma|_\Omega$ such that $\sigma$ has exactly one critical
point $p_j\in\Omega$ with $t_{j-1}<\rho(p_j)<t_j$ for each $j=1,2,\ldots,m$.
Let $k_j$ denote the Morse index of $\sigma$ at $p_j$;
thus $k_j\le n$ for all $j$. By Morse theory \cite{Mi}
the domain $D_{t_j}$ is diffeomorphic to a smooth handlebody obtained
by attaching a handle of index $k_j$ to $D_{t_{j-1}}$ and smoothing
the corners (fig.\ 3). 

%
%
%
%
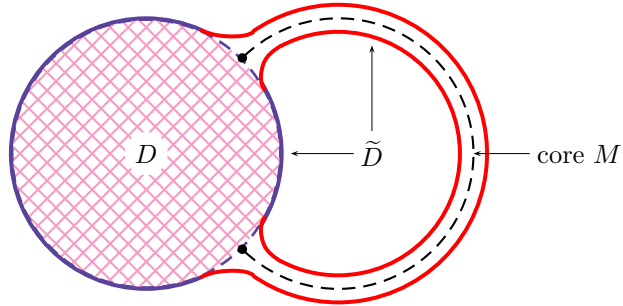
\begin{figure}[ht]
\psset{unit=0.6cm} 
\begin{pspicture}(-4,-4)(9,4)

\pscircle[linestyle=dashed,linewidth=1pt,linecolor=Violet,fillstyle=crosshatch,hatchcolor=Lavender](0,0){3}       
\psarc[linewidth=1.5pt,linecolor=Violet](0,0){3}{66}{294}                                   
\psarc[linewidth=1.5pt,linecolor=Violet](0,0){3}{-32}{32}
\pscircle[linestyle=none,fillstyle=solid,fillcolor=white](0,0){0.5}
\rput(0,0){$D$}

\psarc[linestyle=dashed,arrows=*-*](4.25,0){3}{-135}{135}                        

\psarc[linewidth=1.5pt,linecolor=red](4.25,0){2.7}{-118}{118}                       
\psarc[linewidth=1.5pt,linecolor=red](4.25,0){3.3}{-125}{125}                       

\psecurve[linewidth=1.5pt,linecolor=red](6,1)(2.64,1.4)(2.6,2)(3,2.4)(4,2.7)        
\psecurve[linewidth=1.5pt,linecolor=red](6,-1)(2.64,-1.4)(2.6,-2)(3,-2.4)(4,-2.7)   
\psecurve[linewidth=1.5pt,linecolor=red](0,4)(1.2,2.75)(2,2.6)(2.38,2.7)(4,5)       
\psecurve[linewidth=1.5pt,linecolor=red](0,-4)(1.2,-2.75)(2,-2.6)(2.38,-2.7)(4,-5)  

\psline[linewidth=0.3pt]{<-}(3.2,0)(4.6,0)
\psline[linewidth=0.3pt]{->}(5,0.5)(5,2.4)
\psline[linewidth=0.3pt]{<-}(7.25,0)(8.5,0)

\rput(9.6,0){core $M$}

\rput(5,0){$\wt D$}

\end{pspicture}
\caption{A handlebody $\wt D$.} 
\label{Fig3}
\end{figure}

Recall that a $k$-handle attached to a compact smoothly bounded domain 
$D\subset X$ is a diffeomorphic image of  
$\Delta_{k}\times \Delta_{2n-k} \subset \R^{k}\times\R^{2n-k}$, 
where $\Delta_k$ denotes the closed unit ball in $\R^k$.
The set $b\Delta_k \times \Delta_{2n-k} = S^{k-1}\times \Delta_{2n-k}$ 
gets attached to $bD$, the image of $\Delta_k \times \{0\}^{2n-k}$ 
is called the {\em core disc} (or simply the core) of the handle, 
and the union of $D$ with the handle, suitably smoothed at the corners,
is a {\em handlebody} $\wt D$ shown on fig.\ 3. 
(In practice one often glues a handle to a thickening of $D$.) 
The Morse theory \cite{Mi} tells us that every smooth manifold
is obtained by successive gluing of handles, i.e., it admits 
a {\em handlebody decomposition}.

We are now ready to complete the proof of lemma \ref{Mainlemma}.
For consistency of notation set $W_0:=D_0$ and $J'_0=J$. 
By what has been said, $D_{t_1}$  is diffeomorphic to a handlebody 
$W_1\subset D_{t_1}$ obtained by attaching to $W_0$ a handle of index $k_1$. 
Since $W_0$ is strongly $J$-pseudoconvex and $k_1\le n\ne 2$,
Eliashberg's results from \cite{E} show that the core disc $M$ 
of the handle can be chosen $J$-totally real in $X$ and such that 
its boundary sphere $bM$ is a $J$-Legendrian (complex tangential) 
submanifold of $bW_0$. (See lemma 3.1 in \cite{FS} for 
details of this construction. It is here that the hypothesis $\dim_\C X\ne 2$
is needed; in the exceptional case $\dim_\C X\ne 2$ and $k_1=2$ 
it is in general impossible to find an {\em embedded} totally real 
core disc $M$ for the $2$-handle as is shown by the gauge theory; 
see \cite{E} and \cite{Go1}. We shall discuss this in \S 4 below.) 

After a small homotopic deformation of the almost complex structure 
$J=J'_0$ supported in a neighborhood of the core disc $M$ 
(and away from the set $W_0$) we find a new almost complex 
structure $J'_1$ on $X$ which is integrable near $W_0\cup M$ and
agrees with $J'_0$ near $A\cup W_0$, 
and the handlebody $W_1$ (a thickening of $W_0\cup M$)
can be chosen such that $bW_1$ is smooth strongly $J'_1$-pseudoconvex, 
${\rm Int} W_1$ is $J'_1$-Stein, and $W_0$ is $J'_1$-holomorphically convex 
in $W_1$. If $J$ is integrable on $X$ 
then the same can be accomplished without a homotopic correction 
of $J$ by choosing the core disc $M$ to be real analytic (see \cite{E} and \cite{FS}).

In addition, lemma 5.1 in \cite{FS} shows that we can choose 
$W_1$ sufficiently thin around $W_0\cup M$ such that there exists a map 
$g_1 \colon X\to Y$ which is $J'_1$-holomorphic in a neighborhood of $W_1$ 
and satisfies the following properties:
\begin{itemize}
\item[(a)]  $\sup_{x\in W_0} d_Y\bigl(f(x),g_1(x)\bigr) < \frac{\epsilon}{m}$,
\item[(b)]  $g_1|_A = f|_A$,  
\item[(c)]  $g_1$ is homotopic to $f$ by a homotopy $\{g_t\}_{t\in[0,1]}$ 
consisting of maps defined near $W_1$ which agree with $f$ on $A$, they are 
$J$-holomorphic in a neighborhood of $W_0$, and each of them is $\frac{\epsilon}{m}$-close 
to $f$ on $W_0$. 
\end{itemize}

To obtain the interpolation conditions in (b) and (c)
which are not explicitly stated by lemma 5.1 in \cite{FS}, the reader should
observe that the proof of that lemma relies on theorem 3.2 in 
\cite[p.\ 1924]{FFourier1} which includes interpolation 
on a complex subvariety.

Using the homotopy $\{g_t\}$ we can patch all these maps with $f$ outside 
a certain neighborhood of $W_0$ in order to get a homotopy of
global maps $X\to Y$.

We now proceed to the next set $D_{t_2}$. By the same argument as above,
$D_{t_2}$ is diffeomorphic to a handlebody obtained from $D_{t_1}$ by attaching 
a handle of index $k_2$. As $D_{t_1}$ is diffeomorphic to $W_1$, 
$D_{t_2}$ is also diffeomorphic to a handlebody $W_2 \subset D_{t_2}$ 
obtained by attaching a handle of index $k_2$ to $W_1$. By repeating the 
above arguments we can modify $J'_1$ near the core disc $M_1$ of the handle 
to a structure $J'_2$ which is integrable near $W_1\cup M_1$,
and we then choose $W_2$ to be strongly $J'_2$-pseudoconvex, 
with  $J'_2$-Stein interior, and such that $W_1$ is $J'_2$-holomorphically convex in $W_2$.
After shrinking $W_2$ around $W_1\cup M_1$ we also get a map 
$g_2\colon X\to Y$ which is $J'_2$-holomorphic in a neighborhood 
of $W_2$, it agrees with $f$ on $A$, it satisfies 
$\sup_{x\in W_1} d_Y\bigl(g_2(x),g_1(x)\bigr) < \frac{\epsilon}{m}$,
and it is homotopic to $g_1$ by a homotopy $\{g_t\}_{t\in [1,2]}$ which is fixed on $A$
such that each $g_t$ is $J'_1$-holomorphic near $W_1$ and is 
uniformly $\frac{\epsilon}{m}$-close to $g_1$ on $W_1$. 

Continuing inductively we obtain after $m$ steps 
a handlebody $W_m \subset L$, diffeomorphic to $L$,
and an almost complex structure $J'_m$ on $X$ which 
is integrable in a neighborhood of $A\cup W_m$ and which 
agrees with $J=J'_0$ in a neighborhood of $A\cup  W_0$ 
(in fact, the two structures are homotopic by a homotopy 
that is fixed near $A\cup W_0$), such that $W_m$ is strongly 
$J'_m$-pseudoconvex and its interior is $J'_m$-Stein. 
We also obtain a map $g_m \colon X\to Y$ which is $J'_m$-holomorphic in 
a neighborhood of $W_m$, it agrees with $f$ on $A$, and it satisfies
$\sup_{x\in D_0} d_Y\bigl(f(x),g_m(x)\bigr) < \epsilon$.
The construction also gives a homotopy of maps $X\to Y$ 
connecting $f$ to $g_m$ such that the homotopy is fixed on $A$, 
each map in the family is $J$-holomorphic 
in a neighborhood of $D_0=W_0$ and is uniformly 
$\e$-close to $f$ on $D_0$ (and hence on $K$). 

Our construction of the handlebodies $W_1,\ldots,W_m$ insures 
that there is a diffeomorphism $h\colon X\to X$ such that $h(L)=W_m$
and $h$ is diffeotopic to $id_X$ by a diffeotopy 
that is fixed in an open neighborhood of $A\cup K$.
(We may even insure that $h(D_{t_j})=W_j$ for $j=0,1,\ldots,m$.)

Set $J_1=h^*(J'_m)$ and $f_1=g_m\circ h\colon X\to Y$.
The definition of $J_1$ is equivalent to $dh\circ J_1= J'_m \circ dh$,
which means that $h\colon (X,J_1)\to (X,J'_m)$ is a biholomorphic map.
Note that  $J_1$ is integrable
in a neighborhood of $A\cup L$ (since $J'_m$ is integrable
near $A\cup W_m$), and $J_1$ coincides with $J$ 
near $A$ (since $h$ is the identity near $A$). 

If $\{h_t\}_{t\in[0,1]}$ is a diffeotopy on $X$ 
from $h_0=id_X$ to $h_1=h$ which is fixed near $A\cup K$ 
then $J_t=h_t^*(J'_m)$ is a homotopy
of almost complex structures which is fixed in a neighborhood 
of $A\cup K$ and which connects $J_0=J$ to $J_1$.

If $J$ is integrable on $X$ then $J=J'_0=J'_1=\cdots J'_m$
by the construction, and hence the structure $J_t$ 
is integrable for every $t\in[0,1]$
since conjugation by a diffeomorphism preserves integrability.
This verifies properties (i) and (ii) in lemma \ref{Mainlemma}.

The set $L=h^{-1}(W_m)$ is strongly 
$J_1$-pseudoconvex and its interior is $J_1$-Stein  
since $W_m$ enjoys these properties with respect to $J'_m$. 
Also, $W_j$ is $J'_{j+1}$-holomor\-phically convex in 
$W_{j+1}$ and $J'_j=J'_{j+1}$ near $W_j$ for $j=0,1,\ldots,m-1$;
since $K$ is $J'_0$-holomorphically convex in $U$ and hence in $W_0$, 
we see that $K$ is $J'_m$-holomorphically convex in $W_m$.
Thus $K$ is $J_1$-holomorphically convex in $L$ and hence (iii) holds. 

The map $f_1= g_m\circ h \colon X\to Y$ 
is $J_1$-holomorphic near $L$ 
(since $h\colon (X,J_1)\to (X,J'_m)$ is biholomorphic
and $g_m$ is $J'_m$-holo\-mor\-phic in a neighborhood of $W_m=h(L)$),
so (vi) holds. By the construction we also have
$\sup_{x\in K} d_Y\bigl(f(x),f_1(x)\bigr) < \epsilon$.
A homotopy from $f=f_0$ to $f_1$ satisfying properties (iv) 
and (v) is obtained by combining the individual homotopies
obtained in the construction. This completes the proof.
\end{proof}

\begin{remark} 
H.\ Hamm proved in \cite{Ha1} and \cite{Ha2} that for every 
$n$-dimensional Stein space $X$ and closed complex subvariety
$A\subset X$ the pair $(X,A)$ is homotopically equivalent 
to a relative CW complex of dimension $\le n=\dim_\C X$. 
(The absolute version with $A=\emptyset$ is a well known 
theorem of Lefshetz \cite{Lef}, Abraham and Fraenkel \cite{AF} 
and Milnor \cite{Mi}.) In his proof Hamm used Morse theory for 
manifolds with boundary. The essential step is the following
\cite[pp.\ 2--5]{Ha2}:

{\em Assume that $A$ is a closed complex subvariety of an 
$n$-dimensional Stein space $X$ such that $X\bs A$ is regular 
(without singularities). Let $K\subset L$ be sublevel sets 
of a real analytic, strongly plurisubharmonic Morse exhaustion function 
on $X$. Then $(A\cap L)\cup K$ admits a thickening $D\subset L$ such that 
$A\cup L$ is obtained from $A\cup D$ by attaching handles of index $\le n$.} 

The geometric device in the proof of our lemma \ref{Mainlemma},
using the family of domains $\{D_t\}_{t\in[0,1]}$ 
which increase from $D_0$ to $D_1=L$, 
accomplishes this by using only the classical Morse theory 
for manifolds without boundary.
\end{remark}

%
%
%
%
\section{Proof of theorems \ref{Tmain1} and \ref{Tmain2}}
Theorem \ref{Tmain1} corresponds to the special case 
$K=\emptyset$ and $J=J_X$ of the following.

\begin{theorem}
\label{T3.1}
Let $(X,J_X)$ be a Stein manifold with $\dim_\C X\ne 2$,
let $K\subset X$ be a compact $\cH(X,J_X)$-convex subset 
with smooth strongly $J_X$-pseudoconvex boundary, and  
let $A$ be a closed complex subvariety of $X$.
Assume that $J$ is an almost complex structure on $X$
which is integrable in an open neighborhood of $A\cup K$, 
it agrees with $J_X$ in a neighborhood of $A$, 
and such that $K$ is a strongly $J$-pseudoconvex with 
$J$-Stein interior. Let $Y$ be a complex manifold with 
a distance function $d_Y$ induced by a Riemannian metric.

Given a continuous map $f\colon X\to Y$ which is $J$-holomorphic 
in a neighborhood of $K$ and such that $f|_A\colon A\to Y$ 
is holomorphic, there exists for every $\e>0$ 
a homotopy of pairs $(J_t,f_t)$ $(t\in[0,1])$, where $J_t$ is
an almost complex structure on $X$ and $f_t\colon X\to Y$
is a continuous map, satisfying the following:
\begin{itemize}
\item[(i)] $J_0=J$, and $J_t$ agrees with $J$ 
in a neighborhood of $A\cup K$ for every $t\in [0,1]$,
\item[(ii)]  the structure $J_1$ is integrable Stein on $X$
and $K$ is $\cH(X,J_1)$-convex,
\item[(iii)] $f_0=f$, and $f_t|_A=f|_A$ for every $t\in[0,1]$,
\item[(iv)]  for each $t\in [0,1]$ the map $f_t$ is 
$J$-holomorphic in a neighborhood of $K$ and satisfies
$\sup_{x\in K} d_Y\bigl(f_t(x),f(x)\bigr) <\e$, and
\item[(v)]   the map $f_1\colon X\to Y$ is $J_1$-holomorphic.
\end{itemize}
If $J$ is integrable on $X$ then $J_t$ can be chosen integrable
for every $t\in[0,1]$.
\end{theorem}

We emphasize that the almost complex structure $J$ need not
be homotopic to $J_X$. In fact, the Stein structure $J_X$ is only 
used to obtain a correct handlebody decomposition of the pair 
$(X,A)$ (see remark \ref{JX} below).

\begin{proof}
Since $K$ is strongly $J_X$-pseudoconvex and $\cH(X,J_X)$-convex,
there exists a smooth strongly $J_X$-plurisubharmonic Morse exhaustion function 
$\rho\colon X\to \R$ such that $K=\{x\in X\colon \rho(x)\le 0\}$
and $d\rho\ne 0$ on $bK=\{\rho=0\}$. 
Choose  a sequence $c_0=0<c_1<c_2\ldots$ consisting of 
regular values of $\rho$, with $\lim_{j\to\infty} c_j = +\infty$.
Let $K_j=\{x\in X\colon \rho(x)\le c_j\}$. Set $f_0=f$ and $J_0=J$.
Applying lemma \ref{Mainlemma} we can  inductively construct 
sequences of maps $f_j\colon X\to Y$ and of almost complex
structures  $J_j$ satisfying the following for $j=1,2,\ldots$: 
\begin{itemize}
\item[(a)] $J_j$ is integrable in a neighborhood of $A\cup K_j$
and it agrees with $J_{j-1}$ in a neighborhood of $A\cup K_{j-1}$,
\item[(b)] $K_j$ is strongly $J_j$-pseudoconvex with $J_j$-Stein interior,
and $K_{j-1}$ is $J_j$-ho\-lo\-mor\-phi\-cally convex in $K_j$,
\item[(c)] there is a homotopy of almost complex structures 
$J_{j,s}$ $(s\in[0,1])$, with $J_{j,0}=J_{j-1}$ and $J_{j,1}=J_j$,
which is fixed in a neighborhood of $A\cup K_{j-1}$,
\item[(d)] the map $f_j\colon X\to Y$ is $J_j$-holomorphic
in a neighborhood of $K_j$  and $f_j|_A=f|_A$, and
\item[(e)] there is a homotopy $f_{j,s}\colon X\to Y$ $(s\in[0,1])$
which is fixed on $A$ such that $f_{j,0}=f_{j-1}$, $f_{j,1}=f_j$,
and for every $s\in [0,1]$ the map  $f_{j,s}$ is $J_{j-1}$-holomorphic 
in a neighborhood of $K_{j-1}$ and it satisfies
$$
	\sup_{x\in K_{j-1}} d_Y\bigl(f_{j,s}(x),f_{j-1}(x)\bigr) < 2^{-j-1}\epsilon.
$$
\end{itemize}

Indeed, assuming that we have already constructed the above sequences up to $j-1$,
it suffices to apply lemma \ref{Mainlemma} with $K=K_{j-1}$,
$L=K_j$, $f=f_{j-1}$, $J=J_{j-1}$, and $\e$ replaced by $2^{-j-1}\e$
to get the next complex structure $J_j$ and the next map
$f_j$ satisfying the stated properties.

Condition (a) insures that the limit $\wt J=\lim_{j\to\infty} J_j$
exists and is an integrable complex structure on $X$ which agrees with 
$J$ in a neighborhood of $A\cup K$. The manifold $X$ is exhausted 
by the sequence of strongly $\wt J$-pseudoconvex domains $K_j$ 
with $\wt J$-Stein interior. Property (b) implies that $K_j$ 
is $\cH(X,\wt J)$-convex for $j=0,1,2,\ldots$ 
and hence the manifold $(X,\wt J)$ is Stein. By combining the 
individual homotopies furnished by (c) we obtain a homotopy of 
almost complex structures on $X$ which connects $J$ to $\wt J$
and which is fixed in a neighborhood of $A\cup K$.

Properties (d) and (e) insure that the sequence of maps $f_j\colon X\to Y$ 
converges uniformly on compacts in $X$ to a $\wt J$-holomorphic 
map $\wt f=\lim_{j\to \infty} f_j \colon X\to Y$ satisfying
$\wt f|_A=f|_A$ and $\sup_{x\in K} d_Y\bigl(\wt f(x),f(x)\bigr) <\e$.
Finally, condition (e) implies that the homotopies
$f_{j,s}$ $(s\in[0,1], j=1,2,\ldots)$ can be assembled
into a homotopy from $f$ to $\wt f$ which is fixed on $A$, 
holomorphic on $K$, and $\e$-close to $f$ on $K$. 

Changing the notation so that $\wt J$ is denoted $J_1$
and $\wt f$ is denoted $f_1$ we obtain the conclusion of 
theorem \ref{T3.1}.
\end{proof}

\begin{remark}
\label{JX}
The Stein structure $J_X$  was used in the above proof 
only to insure that for every $j=1,2,\ldots$  there is a thickening 
$D_{j-1} \subset K_j$ of the set $K_{j-1}\cup (A\cap K_j)$ 
such that $A\cup K_j$ is obtained (up to a diffeomorphism) 
by attaching handles of index $\le \dim_\C X$ to $A\cup D_{j-1}$. 
(In the proof of lemma \ref{Mainlemma} this  was
shown using the notation $K_j=L$, $K_{j-1}=K$ and $D_{j-1}=D_0$.) 
This leads to a proof of theorem \ref{Tmain1} under the weaker conditions 
that $(X,J)$ is an almost complex manifold of real dimension $2n\ne 4$ 
such that $J$ is integrable in a neighborhood of a closed 
Stein subvariety $A\subset X$, and $X$ is exhausted 
by an increasing sequence  of compact strongly $J$-pseudoconvex 
domains $K_0\subset K_1\subset \ldots\subset \cup_{j=0}^\infty K_j =X$
such that every pair $(A\cup K_j,A\cup K_{j-1})$ satisfies 
the above topological condition. 
\end{remark}

{\em Proof of theorem  \ref{Tmain2}.} 
We shall use the same tools as in the proof of theorem \ref{T3.1},
but will change the induction procedure. Unlike in theorem \ref{T3.1},
the complex structure on $X$ will remain fixed during the 
entire  proof.

Let $K_0\subset K_1\subset \cdots \subset \cup_{j=0}^\infty K_j = X$ 
be an exhaustion of $X$ by compact, smoothly bounded,
strongly pseudoconvex sets as in the proof of theorem \ref{T3.1}. 
Set $f_0=f$. We shall assume that $f_0$ is holomorphic in a neighborhood 
of $K_0$ (choosing $K_0=\emptyset$ if so desired.) 
Let $d_Y$ be a distance function on $Y$.

Given an $\e>0$ we shall inductively construct a sequence of compact, 
smoothly bounded, strongly pseudoconvex sets 
$\emptyset=O_{-1} \subset O_0\subset O_1\subset \ldots \subset X$,
a sequence of smooth diffeomorphisms $h_j\colon X\to X$, 
and a sequence of maps $f_j\colon X\to Y$ satisfying the 
following properties for  $j=1,2,\ldots$:
\begin{itemize}
\item[(i)]    $h_j(K_j)=O_j$, and $h_j$ is diffeotopic to $h_{j-1}$ by a diffeotopy 
which is fixed in a neighborhood of $A\cup K_{j-1}$,
\item[(ii)]   $O_{j-1}$ is holomorphically convex in $O_{j}$,
\item[(iii)]   $f_j$ is holomorphic in an open neighborhood of $O_j$
and satisfies $f_j|_A=f|_A$,
\item[(iv)]    there is a homotopy $f_{j,s}\colon X\to Y$ $(s\in [0,1])$ 
such that $f_{j,0}=f_{j-1}$, $f_{j,1}=f_j$, the homotopy is fixed on $A$,
each map $f_{j,s}$ $(s\in[0,1])$ is holomorphic in a neighborhood of $O_{j-1}$, and 
$$
	\sup_{x\in O_{j-1}} d_Y\bigl(f_{j,s}(x), f_{j-1}(x)\bigr) < 2^{-j-1}\e,
	\quad s\in[0,1].
$$
\end{itemize}

We begin by setting $O_0=K_0$, $h_0=id_X$ and $f_{0,s}=f_0$ for all $s\in[0,1]$.
Suppose inductively that we have already constructed our sequences 
up to an index $j\in \Z_+$; thus the map $f_j\colon X\to Y$ is holomorphic
on $A$ and in an open neighborhood of $O_j$. Property (i) implies that
$h_j$ equals the identity map in a neighborhood of $A\cup K_0$.
Hence $O_j\cap A=K_j\cap A$, and this set is holomorphically convex in $A$
since $K_j$ is $\cH(X)$-convex. The set $O_j$, being strongly pseudoconvex,
admits a basis of open Stein (strongly pseudoconvex) neighborhoods in $X$.
In this situation theorem 3.1 in \cite{FFourier2} applies and furnishes 
a  map $f'_j\colon X\to Y$ which is holomorphic in an open neighborhood
$V_j\supset A\cup O_j$ and which approximates $f_j$ as close as desired
uniformly on $O_j$. Replacing $f_j$ by $f'_j$ we may therefore
assume that $f_j$ is holomorphic in an open set $V_j\supset A\cup O_j$.

Applying lemma \ref{Mainlemma} with $f=f_j$, $K=K_j$ and $L=K_{j+1}$
we find a compact domain $D_j \subset K_{j+1}$ with strongly pseudoconvex 
boundary (denoted $D_0$ in lemma \ref{Mainlemma})
such that $(A\cap K_{j+1})\cup K_j \subset D_j$, $K_{j+1}$ is obtained 
from $D_j$ by attaching finitely many handles of index $\le n=\dim_\C X$,
and $h_j(D_j) \subset V_j$. The last inclusion is trivially satisfied
in a neighborhood of $A$ where $h_j$ coincides with the identity map,
while outside this neighborhood $D_j$ can be chosen as close as desired
to $K_j$; since  $h_j(K_j)=O_j \subset V_j$, the inclusion follows.

Set $O'_j=h_j(D_j)$. If the above approximations were 
chosen sufficiently close then $O'_j$ is a compact set with smooth 
strongly pseudoconvex boundary (since $bO'_j$ coincides with $bD_j$
near the subvariety $A$, and elsewhere $bO'_j$ is $\cC^\infty$-close 
to the strongly pseudoconvex hypersurface $h_j(b K_j) = bO_j$).  
Note that $O_j$ is holomorphically convex in $O'_j$ provided that 
$D_j$ is chosen in a sufficiently small neighborhood of $(A\cap K_{j+1})\cup K_j$.
Applying the diffeomorphism $h_j$ to the above sets we see that
$h_j(K_{j+1})$ is diffeomorphic to a handlebody $O_{j+1}$ obtained 
from $O'_j=h_j(D_j)$ by attaching finitely many handles of index $\le n$. 

We now proceed as in the proof of theorem \ref{T3.1}.
By Lemma 5.1 in \cite{FS} the above handles can be chosen such 
that the resulting handlebody $O_{j+1}$ has smooth strongly pseudoconvex boundary,
$O'_j$ is holomorphically convex in $O_{j+1}$,  and there is a map
$f_{j+1}\colon X\to Y$ which is holomorphic in a neighborhood of 
$O_{j+1}$, it agrees with $f_j$ on $A$, and 
$\sup_{x\in O_{j}} d_Y\bigl(f_{j+1}(x), f_{j}(x)\bigr) < 2^{-j-2}\e$. 
The same lemma provides a homotopy from $f_j$ to $f_{j+1}$ satisfying 
property (iv) for the index $j+1$.

Since $O_{j+1}$ is constructed from $O'_j$ by using the topological data 
provided by the pair $D_j\subset K_{j+1}$ and since all handles used in 
the construction of $O_{j+1}$ are contained in $X\bs A$, 
there exists a diffeomorphism 
$g_j\colon X \to X$ which maps $h_j(K_{j+1})$ onto $O_{j+1}$ and 
which is diffeotopic to $id_X$ by a diffeotopy which is 
fixed (equal the identity map) in a neighborhood of $A\cup O'_j$. 
The map $h_{j+1}=g_j\circ h_j \colon X\to X$ is a diffeomorphism 
of $X$ which maps $K_{j+1}$ onto $O_{j+1}$ and is diffeotopic to $h_j$ 
by a diffeotopy which is fixed near $A\cup K_j$. 
The induction may now continue.

Properties (i)--(iv) insure that $\Omega=\cup_{j=0}^\infty O_j \subset X$
is a Stein domain which contains $A\cup K_0$, and the sequence
$f_j$ converges uniformly on compacts in $\Omega$ to a holomorphic map
$\wt f = \lim_{j\to\infty}f_j \colon \Omega \to Y$
satisfying $\wt f|_A=f|_A$ and $\sup_{x\in K_0} d_Y\bigl( \wt f(x),f(x)\bigr) <\e$.
Also, there is a homotopy of maps $\Omega\to Y$ from $f|_\Omega$ to $\wt f$ 
which is holomorphic on $K_0$ and is $\e$-close to $f_0$ on $K_0$.
Property (i) also gives a diffeomorphism 
$h = \lim_{j\to\infty} h_j\colon X\to h(X) = \Omega$
which is diffeotopic to $id_X$ and which equals
the identity map in a neighborhood of $A$. 
It follows that the map $\wt f\circ h\colon X\to Y$ is homotopic 
to $f$, thereby  completing the proof of theorem \ref{Tmain2}.

\section{The case $\dim_\C X=2$}
The proof of lemma \ref{Mainlemma}, and hence of theorems \ref{Tmain1} and \ref{Tmain2},
breaks down when $X$ is a Stein surface ($\dim_\C X=2$), the reason being that 
a certain framing obstruction may arise when trying to add a 2-handle 
with an embedded totally real core disc attached along a Legendrian knot
to a given strongly pseudoconvex boundary in $X$. 
This obstruction in the proof has been pointed out by Eliashberg \cite{E}, 
and the Seiberg-Witten theory subsequently confirmed 
that it cannot be removed in general.
In particular, there exist smooth, orientable, almost complex $4$-manifolds 
$(X,J)$ with a handlebody decomposition without handles of index $>2$ 
which  do not admit any Stein structure; one such example
is the manifold $X=S^2\times\R^2= \C\P^1\times \C$. (Many futher 
examples can be found in \cite{Go1}.)
The key obstruction for the existence of a Stein structure is provided by the 
{\em generalized adjunction inequality} which states that for every
closed, orientable, smoothly embedded 2-surface $S$ in a Stein
manifold $X$, with the only exception of a null-homologous 2-sphere, we have
$$
	[S]^2 + |c_1(X)\cdotp S| \le -\chi(S).
$$
(See Chapter 11 in \cite{GS}, or \cite{OS}, for a proof, references to the
original papers and further results.)
Conversely, a closed embedded orientable 2-surface in an arbitrary
complex manifold $X$ which satisfies the above inequality 
is isotopic to another embedding with a basis of 
tubular open Stein neighorhoods in $X$ \cite{Steinnbds}.

On the other hand, Gompf proved that there always exist {\em exotic Stein structures}
on any such 4-manifold $X$ \cite{Go1}, \cite{Go2}. More precisely, given 
a smooth, almost complex 4-manifold $(X,J)$ with a Morse exhaustion function
without critical points of Morse index $>2$, there exist a Stein surface 
$(X',J')$ and an orientation preserving {\em homeomorphism} $h\colon X \to X'$ 
such that the class determined by the almost complex structure $J'$ via $h$
agrees with the class of $J$. 

Keeping the same hypotheses on $(X,J)$, the authors have shown in \cite[\S 7]{FS} 
that for any continuous map $f\colon X\to Y$ to a complex manifold $Y$, 
a Stein surface $(X',J')$ and a homeomorphism $h\colon X\to X'$
in Gompf's theorem can be chosen such that there exists a $J'$-holomorphic 
map $f'\colon X'\to Y$ with the property that the map 
$\wt f= f'\circ h\colon X\to Y$  is homotopic to $f$.
If in addition the almost complex structure $J$ on $X$ is integrable 
(but not necessarily Stein), one can realize such $(X',J')$ as an open $J$-Stein domain 
$\Omega \subset X$ which is homeomorphic to $X$ (theorem 1.2 in \cite{FS}; 
without considering  mappings this is again due to Gompf \cite{Go2}). 

The constructions in \cite{Go1}, \cite{Go2} and \cite{FS} use 
{\em kinky discs} and {\em Casson handles} at every place where a 
framing obstruction arises in the construction, together with
the famous result of Freedman to the effect that a Casson handle
is homeomorphic to a standard index two handle $\Delta_2\times \Delta_2\subset\R^4$
\cite{Freedman}, \cite{FQ}. By using the same tools, together with the methods
explained in this paper, one can prove the following interpolation
theorem which is the analogue of theorem \ref{Tmain2} in the
case $\dim_\C X=2$.

\begin{theorem}
\label{T4.1}
Let $X$ be a Stein surface and let $A$ be a closed complex subvariety of $X$.
Given a continuous map $f\colon X\to Y$ to a complex manifold $Y$ 
such that $f|_A\colon A\to Y$ is holomorphic, there exist a 
Stein domain $\Omega\subset X$ containing $A$, a holomorphic map 
$\wt f\colon \Omega\to Y$, and an orientation preserving homeomorphism 
$h\colon X\to \Omega$ which is homeotopic to $id_X$ by a homeotopy
that is fixed on a neighborhood $A$, such that the map $\wt f\circ h\colon X\to Y$
is homotopic to $f$ relative to $A$.
\end{theorem}

This can be proved by modifying the proof of theorem \ref{Tmain2} 
in \S 3 above, and the necessary modification is explained in 
the proof of theorem 1.2 in \cite[\S 7]{FS}.  To avoid excessive 
repetition we shall only recall the essential points.

Let $J$ denote the given Stein structure on $X$.
We assume the notation used in the proof of theorem \ref{Tmain2} 
in \S 3 above. In that proof it is explained  how one obtains
a strongly pseudoconvex handlebody $O_{j+1}$ by attaching 
handles of index $\le n$ to  a strongly pseudoconvex domain $O'_j$.
Each of the handles must have an embedded totally real core disc
whose boundary circle is attached to 
the previous strongly pseudoconvex hypersurface along
a Legendrian knot; this enables us to choose the next handlebody
to be strongly pseudoconvex and to approximate the holomorphic map
by a map which is holomorphic on a neighborhood of the new (larger)
handlebody.

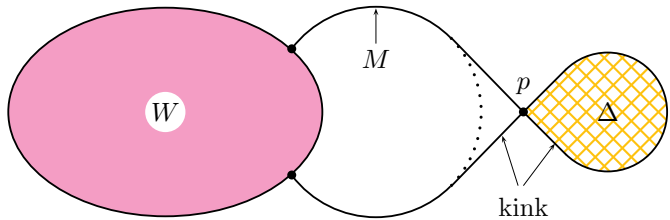
\begin{figure}[ht]
\psset{unit=0.7cm, linewidth=0.7pt} 

\begin{pspicture}(0,-3)(16,3)
\definecolor{myblue}{rgb}{0.66,0.78,1.00}

%
%
\psellipse[fillstyle=solid,fillcolor=Lavender](4,0)(3,2)
\pscircle[fillstyle=solid,fillcolor=white,linestyle=none](4,0){0.4}
\rput(4,0){$W$}

%
%
\psarc(8,0){2}{45}{143}
\psarc(8,0){2}{-143}{-45}
\psarc[linestyle=dotted,linewidth=1.2pt](8,0){2}{-45}{45}
\psline(9.42,1.42)(10.8,0)   
\psline(9.42,-1.42)(10.8,0)

%
%
\pscustom[fillstyle=crosshatch,hatchcolor=Apricot]
{
\psline(10.8,0)(11.6,-0.8)
\psarc(12.4,0){1.13}{-135}{135}
\psline(11.6,0.8)(10.8,0)
}

\psdots(10.8,0)(6.4,1.2)(6.4,-1.2)

\rput(8,1){$M$}
\psline[linewidth=0.2pt]{->}(8,1.3)(8,1.97)
\rput(12.4,0){$\Delta$}

\psline[linewidth=0.2pt]{<-}(10.4,-0.45)(10.7,-1.4)
\psline[linewidth=0.2pt]{->}(10.8,-1.4)(11.4,-0.6)
\rput(10.8,-1.8){kink}
\rput(10.8,0.5){$p$}

\end{pspicture}
\caption{A kinky disc $M$ with a trivializing 2-cell $\Delta$}
\label{Fig4}
\end{figure}

When $\dim_\C X=2$, a framing problem may arise for handles of index 2,
and a required totally real embedded core disc $M$ does not exist in general.
As explained in \cite{FS} (and before that in \cite{Go1}), the
problem can be resolved by choosing an embedded core disc $M$ 
which is attached to the given strongly pseudoconvex domain $W\subset X$
along a Legendrian knot $bM\subset bW$, and then adding finitely many 
{\em positive kinks} to $M$. More precisely, we remove 
from $M$ finitely many small pairwise disjoint discs and 
glue along each of the resulting circles an immersed disc 
with one positive double point. (Fig.\ \ref{Fig4}, borrowed from \cite{FS},
shows a kink with a trivializing disc $\Delta$ which 
will be attached at the next step in order to cancel the superfluous loop 
at the double point $p$. A model kink used in \cite{FS} is provided by an explicit
immersed Lagrangian sphere in $\C^2$, due to Weinstein \cite{We}.)  

%
%
%
%

As explained in \cite{FS}, kinking the core disc sufficiently
many times gives an immersed disc which can be deformed 
to a totally real immersed disc $M' \subset X\bs {\rm Int W}$, 
attached to $bW$ along a Legendrian knot $bM' \subset bW$. 
It is then possible to find a thin strongly pseudoconvex 
neighborhood $W'\subset X$ of $W\cup M'$ and a holomorphic map 
$W'\to Y$ which approximates the given
initial map $f\colon X\to Y$ uniformly on $W$.
The manifold $W'$ does not have the correct topology 
(it is not even homeomorphic to the domain obtained by attaching 
to $W$ a standard handle with an embedded core disc). 
The problem is partially corrected in the next stage of the 
construction by attaching to $W'$ a trivializing $2$-disc
$\Delta$ at each of the kinky points in order to cancel
the extra loop. Unfortunately the framing obstruction 
arises at this disc as well, requiring us to place another 
kink on $\Delta$ which will require a new trivializing  disc,
etc. The ensuing procedure is always infinite, it can be carried 
out in a small neighborhood of the initial kinky point in $M$,
and (the main point!) it converges to an attached Casson handle 
which is homeomorphic to the standard 2-handle $\Delta_2\times\Delta_2$
(Freedman \cite{Freedman}, \cite{FQ}). Performing this construction 
inside $X$ gives a Stein domain
$\Omega \subset X$ which is homeomorphic to $X$, but in general not 
diffeomorphic to $X$ due to the presence of Casson handles. 
A more precise description of this construction can be found 
in \cite{FS}, and of course in \cite{Go2} for the topological part. 
To insure that $\Omega$ contains the 
given subvariety $A \subset X$ we follow the proof of 
theorem \ref{Tmain2} with these modifications.

%
%
%
%
%

\bibliographystyle{amsplain}

\end{document}